\documentclass[12pt,pdftex]{amsart}
\pdfoutput=1

\usepackage[margin=2.5cm]{geometry}
\usepackage{hyperref}
\usepackage{amsmath,amsthm,amssymb,enumerate,float,graphicx,multicol,xcolor,tikz}
\usepackage[all]{xy}

\usepackage{tikz}
\usepackage{tikz-3dplot}
\usepackage{pgfplots}
\usepackage{tikz}
\tdplotsetmaincoords{70}{110}

\newcommand{\rotateRPY}[3]
{   \pgfmathsetmacro{\rollangle}{#1}
    \pgfmathsetmacro{\pitchangle}{#2}
    \pgfmathsetmacro{\yawangle}{#3}

    \pgfmathsetmacro{\newxx}{cos(\yawangle)*cos(\pitchangle)}
    \pgfmathsetmacro{\newxy}{sin(\yawangle)*cos(\pitchangle)}
    \pgfmathsetmacro{\newxz}{-sin(\pitchangle)}
    \path (\newxx,\newxy,\newxz);
    \pgfgetlastxy{\nxx}{\nxy};

    \pgfmathsetmacro{\newyx}{cos(\yawangle)*sin(\pitchangle)*sin(\rollangle)-sin(\yawangle)*cos(\rollangle)}
    \pgfmathsetmacro{\newyy}{sin(\yawangle)*sin(\pitchangle)*sin(\rollangle)+ cos(\yawangle)*cos(\rollangle)}
    \pgfmathsetmacro{\newyz}{cos(\pitchangle)*sin(\rollangle)}
    \path (\newyx,\newyy,\newyz);
    \pgfgetlastxy{\nyx}{\nyy};

    \pgfmathsetmacro{\newzx}{cos(\yawangle)*sin(\pitchangle)*cos(\rollangle)+ sin(\yawangle)*sin(\rollangle)}
    \pgfmathsetmacro{\newzy}{sin(\yawangle)*sin(\pitchangle)*cos(\rollangle)-cos(\yawangle)*sin(\rollangle)}
    \pgfmathsetmacro{\newzz}{cos(\pitchangle)*cos(\rollangle)}
    \path (\newzx,\newzy,\newzz);
    \pgfgetlastxy{\nzx}{\nzy};
}
\tikzset{RPY/.style={x={(\nxx,\nxy)},y={(\nyx,\nyy)},z={(\nzx,\nzy)}}}

\newtheorem{thm}{Theorem}

\newtheorem{lem}[thm]{Lemma}

\newtheorem{con}[thm]{Conjecture}

\theoremstyle{definition}
\newtheorem{dfn}[thm]{Definition}

\newtheorem{rem}[thm]{Remark}

\newtheorem{conv}[thm]{Convention}
\newtheorem{nota}[thm]{Notation}

\DeclareMathOperator{\LS}{LS}

\DeclareMathOperator{\Spec}{Spec}

\newcommand{\fX}{\mathfrak{X}}
\newcommand{\fY}{\mathfrak{Y}}
\newcommand{\fU}{\mathfrak{U}}
\newcommand{\fM}{\mathfrak{M}}
\newcommand{\ff}{\mathfrak{f}}

\newcommand{\RR}{\mathbb{R}}
\newcommand{\PP}{\mathbb{P}}
\newcommand{\ZZ}{\mathbb{Z}}
\newcommand{\NN}{\mathbb{N}}
\newcommand{\FF}{\mathbb{F}}

\renewcommand{\AA}{\mathbb{A}}

\newcommand{\shO}{\mathcal{O}}

\newcommand{\mut}{\mu}

\begin{document}

\renewcommand{\thefootnote}{\fnsymbol{footnote}}

\title[Singular Log Structures]{Singular Log Structures and \\ Log
  Crepant Log Resolutions I}

\author{Alessio Corti}
\address{Department of Mathematics, Imperial College London, 180 Queen's Gate, London, SW7 2AZ, UK}
\email{a.corti@imperial.ac.uk}
\thanks{\thanks{A.C.\  received support from EPSRC Programme Grant
  EP/N03189X/1.} \thanks{H.R. is supported by the NFR Fripro grant Shape2030.}}

\author{Tim Graefnitz}
\address{Leibniz-Universit\"at Hannover \\ Institut f\"ur Algebraische Geometrie \\ Welfengarten 1, 30167 Hannover \\ Germany}
\email{graefnitz@math.uni-hannover.de}

\author{Helge Ruddat}
\address{Department of Mathematics and Physics\\ University of Stavanger\\ P.O. Box 8600 Forus\\ 4036 Stavanger\\ Norway}
\email{helge.ruddat@uis.no}


\date{}

\begin{abstract} \noindent We introduce a class of singular log
  schemes in three dimensions and conjecture that log schemes in this
  class admit log crepant log resolutions. We provide examples
  as evidence and relate this conjecture to the conjecture made
  in~\cite{MR4381899} and the Gross--Siebert program.
\end{abstract}


\maketitle

\tableofcontents{}

\section{Introduction}
\label{sec:introduction}

We introduce a class of singular log structures in dimension $3$, which we
call zero-mutable log structures. We conjecture that log structures
in this class admit log
crepant log resolutions.

We view this as a first step in the development of a new subject of \emph{log birational geometry}.
There are several reasons for pursuing this direction; here, we highlight three:

\begin{enumerate}[(1)]
\item From its inception, log geometry cries out for this development. What constitutes a
singular log structure, and what does it mean to resolve it? Can the
methods of the minimal model program be adapted to the context of log
schemes?

\item The class of zero-mutable log structures that is the focus of
  this paper arises naturally when
  attempting to deform a toric Fano \mbox{$3$-fold} to a smooth (or
  terminal) Fano \mbox{$3$-fold}. The starting point is a
  reflexive polytope $Q\subset M$ having the origin in its strict interior, and
  one wants to deform the toric Fano \mbox{$3$-fold} $X_Q$ that has
  $Q$ as its moment polyhedron. A possible way to do this
  is to first degenerate $X_Q$ further: take the central
  subdivision $\star Q$ of $Q$ and degenerate $X_Q$ to the reducible
  Fano \mbox{$3$-fold} $X=X_{\star Q}$ that has $\star Q$ as its
  moment polyhedral complex. Next one endows $X$ with a log structure,
  making it into a log scheme $X^\dagger$. We claim that log
  structures that one wants here for the sake of smoothability are locally zero-mutable. It turns out that
  instead of deforming $X^\dagger$ directly it is easier to deform a log resolution
  $\varepsilon^\dagger \colon Y^\dagger \to X^\dagger$ and that deformations
  of $Y^\dagger$ induce deformations of $X^\dagger$. 
\item Taking this further, we want to compute the
  Gromov--Witten theory of the smoothing from data on $X^\dagger$. Now
  $X^\dagger$ is singular, hence it does not have a well-defined log
  Gromov--Witten theory. One is led to conjecture that $X^\dagger$ has a
  distinguished class of log-resolutions: these will have a
  well-defined log Gromov--Witten theory. From the perspective of intrinsic
  mirror symmetry~\cite{MR4462625}, one ought to construct
  the Landau--Ginzburg mirror dual of $X^\dagger$ from the log Gromov--Witten theory of a crepant resolution of $X^\dagger$.
\end{enumerate}

All varieties and schemes in this paper are defined over a field $k$.

\subsection*{Synopsis}
\label{sec:synopsis}

In \S~\ref{sec:conjecture} we introduce the class of zero-mutable log
structures and state our conjecture. In \S~\ref{sec:examples} we
explain examples supporting the conjecture. In the final
\S~\ref{sec:relation-} we provide context by sketching 
connections with the work of Corti--Filip--Petracci~\cite{MR4381899}
and the Gross--Siebert program~\cite{MR2846484, MR3415066, MR4520304}. The last
section can be read right after \S~\ref{sec:surfaces}, that is,
skipping almost all of \S~\ref{sec:examples}.

\subsection*{Acknowledgements}
\label{sec:acknowledgements}

A large part of this research has been carried out with the help of the hospitality of Mathematisches Forschungsinstitut Oberwolfach.
Some of the ideas in this paper go back to the thesis of Andrea
Petracci at Imperial College London~\cite{Petracci_2017}. We are
grateful to Tom Ducat for helpful suggestions on the explicit
construction of the resolutions in Section~\ref{sec:a_n-singularity}
and for explaining to us some of the work in his Warwick University
thesis \cite{Ducat_2015}.  

\section{The Conjecture}
\label{sec:conjecture}
We introduce the class of zero-mutable log structures by defining their local property and conjecture that log structures of this class permit resolutions that enjoy specific geometric properties.
\subsection{Log data}
\label{sec:log-structures}
A log datum captures the discrete information of a log singularity playing a similar role as the Newton polyhedron for a hypersurface singularity.

\begin{dfn}
  \label{dfn:log datum}
Let $L\cong \ZZ^2$ be a rank-$2$ lattice. 
A \emph{log datum} on $L$ is a finite set
\[
  S=\{(e_i, \nu_i)\mid i\in I \}
\]
of pairs $(e_i, \nu_i)$ where
\begin{enumerate}[(i)]
\item For all $i\in I$, $e_i\in L$, and we write $e_i=\ell_i u_i$ with $u_i$
  primitive and $\ell_i\in \NN_+$. We assume that the $u_i$ are
  pairwise distinct. 
\item For all $i\in I$, $\nu_i\vdash \ell_i$ is a partition. We write
  \[
\nu_i =(\ell_{i,1}\geq \ell_{i,2} \geq \cdots \geq \ell_{i,k(i)}\geq
0) \quad \text{where} \quad \ell_i=\sum_{k=1}^{k(i)} \ell_{i,k}
    \]
\item The datum is subject to the condition:
  $\sum_{i\in I}e_i=0$; equivalently, $\sum_{i\in I} \ell_iu_i=0$.
\end{enumerate}  
\end{dfn}

\begin{dfn}
  \label{dfn:rank-1-log datum}
A log datum with $|I|=2$ we call \emph{rank one} and a log datum with $|I|>2$, we call \emph{rank two}. A rank one log datum takes the form $S=\{(e,\nu_1),(-e,\nu_2)\}$ for a non-zero vector $e\in L$ of index $\ell$ with two partitions $\nu_1,\nu_2\vdash\ell$. We call a rank one log datum \emph{zero-mutable} if the two partitions agree, $\nu_1=\nu_2$.
\end{dfn}

\begin{dfn}
  A log datum on $L$ is \emph{irreducible} if
  \begin{enumerate}[(i)]
  \item $\gcd(\ell_i\mid i \in I)=1$. 
  \item for all $J\subset I$, if $\sum_{i\in J}e_i=0$ then $J=I$. 
  \end{enumerate}
\end{dfn}

For simplicity, in this paper we only work with irreducible log data. 

\begin{rem}
  Condition (iii) in Definition~\ref{dfn:log datum} implies that the vectors $e_i$ in a log datum are the edge vectors of a plane lattice polygon
  $P\subset L$ uniquely determined up to translation. 
It is convenient to illustrate log data and mutations thereof by drawing this polygon, see Figure~\ref{fig:othermut}. 
On the other hand, thinking of log data more like rays in a fan is closer to the description of the toroidal crossing space $X$ on which the log structure will be placed in Notation~\ref{not-tcs} and Lemma~\ref{lem:LogStr} below, see also Figure~\ref{fig:An}. From this latter perspective, Condition (iii) implies that $X$ is Gorenstein.
\end{rem}

\begin{conv}
  \begin{enumerate}[(i)]
  \item Our lattices $L$ will always come with a choice of
    orientation. This orientation endows $L$ with a symplectic form
    \[
\{\cdot, \cdot\} \colon L\times L \to \ZZ
      \]
      that, in every oriented integral basis $\varepsilon_1, \varepsilon_2$ of $L$, is given by the matrix
      \[
        \begin{pmatrix}
          0 & 1 \\
          -1 & 0
        \end{pmatrix}
      \]
      so, for example
      \[
        \{\varepsilon_1,\varepsilon_2\}=-\{\varepsilon_2,
        \varepsilon_1\}=1
        \]
  \item In a log datum $S$, for all $i\in I$,
  we denote by $i+1\in I$ the index such that $e_{i+1}$ immediately follows
  $e_i$ in the anticlockwise direction. 
  \end{enumerate}
 \end{conv}

\begin{nota}
\label{not-tcs}
Let $S$ be a log datum on $L$, let $M=L\oplus \ZZ$, and write
$u=(0,1)\in M$. Consider the
fan $\Sigma$ in $M$ with maximal cones the
cones $\sigma_{i}=\langle u_i, u_{i+1} ,u\rangle_+$, $i\in
I$. We also denote by $\rho_i =\langle u_i, u \rangle_+$
the \emph{walls} of the fan, and by $\omega=\langle u\rangle_+$ the
(unique) \emph{joint}.\footnote{Following the terminology
  of~\cite{MR2213573} and~\cite{2023arXiv231213867C}, a wall is a
  codimension-$1$ cone and a joint is a codimension-$2$ cone}

This fan
is the moment polyhedral complex of an affine \mbox{$3$-fold} that we denote
by $X_S$ or simply $X$. Denoting
by $k[\Sigma]$ the Stanley--Reisner ring of the fan, we have that
$X=\Spec k[\Sigma]$. Note that $X$ is reducible with irreducible
components $X_i=\Spec k[\sigma_i\cap M]$ intersecting
along the surfaces:
\[
D_i=\Spec k[\rho_i\cap M] = X_{i-1} \cap X_i
\]
Finally note that $k[\omega \cap M]=k[u]$, and the surfaces $D_i$ intersect along the curve $\AA^1_u=\Spec
k[u]$, the affine line with coordinate $u$.
\end{nota}

\begin{rem}
  In fact we only care about the germ (in the Zariski topology) of $X$
  at the origin, but we do not bother being precise about this for the sake of focus on the new ideas. 
\end{rem}

We explain how to endow $X$ with a ``reasonable'' log structure over the standard log point $k^\dagger$. We can put a log structure
on $X$ over $k^\dagger$ by the methods of~\cite{MR2213573}: indeed, the components of
$X$ are toric varieties; moreover, in the language of that
paper, $X$ is the ``vanilla gluing'' of its components. Equivalently,
$X$ is naturally a gtc (\emph{g}enerically \emph{t}oroidal \emph{c}rossing)
space~\cite{2023arXiv231213867C} and we can use this paper to endow
$X$ with a log structure that is compatible with the gtc structure (\cite{2023arXiv231213867C}, Definition 4.2).\footnote{The point of the present paper is to
  construct log resolutions. The log resolutions that we consider are
  not union of toric varieties and hence~\cite{2023arXiv231213867C}
 applies but~\cite{MR2213573} does not apply.} Either way, endowing
$X$ with a log structure amounts to the following:

\begin{lem}
  \label{lem:LogStr}
The set $\LS_{k^\dagger}(X)$ of
log structures on $X$ over $k^\dagger$ compatible with the gtc
structure is the set of data consisting of
\begin{enumerate}[(i)]
\item For all $i\in I$, a function $f_i \in k[\rho_i\cap M]$ that is not divisable by $u$ subject
  to the condition:
\item $\prod_{i\in I} (f_i|k[u])^{e_i}=1$.\footnote{We denote by
    $f_i|k[u]$ the image of $f_i$ under the obvious ring homomorphism
    $k[\rho_i \cap M]\to k[u]$. More
    precisely, the formula reads $\prod_{i\in I} (f_i|k[u])\otimes e_i = 1$ in
    $k(u)^\times \otimes L$, a tensor product of $\ZZ$-modules
    (a.k.a.\ abelian groups) the first factor under multiplication, the
    second factor under addition.}   
\end{enumerate}  
\end{lem}

\begin{dfn}
  \label{dfn:wall_fun}
  The functions $f_i\in k[\sigma_i\cap M]$ are called \emph{wall
    functions}. The condition (ii) is called the \emph{joint compatibility condition}.
\end{dfn}

\begin{proof}
  Unpack~\cite[Theorem~3.22]{MR2213573} or~\cite[Theorem~4.3]{2023arXiv231213867C}.
\end{proof}

\begin{conv}
  \label{conv:singular_log_structures}
  \begin{enumerate}[(1)]
  \item  We will allow our wall functions to be regular
    function on $D_i$ and to vanish along divisors (curves)
  \[
    Z_i = (f_i=0) \subset D_i,
  \]
  so $f_i\in\Gamma(D_i,\shO_{D_i})\cap \Gamma(D_i\setminus Z_i,\shO_{D_i\setminus Z_i}^\times)$.
  Denoting $Z=\cup_{i\in I} Z_i$,  we understand the joint compatibility condition to
  hold as a statement about invertible functions on the locus $X\setminus Z$.

  A set of functions $f_i$ satisfying the joint compatibility condition gives a log structure
  $\mathfrak{M}_{X\setminus Z}$ on $X\setminus Z$ over
  $k^\dagger$ and, denoting by $j\colon X\setminus Z \hookrightarrow
  X$ the inclusion, we then take the direct image log structure $j_\star
  \mathfrak{M}_{X\setminus Z}$ on $X$. This log
  structure is badly behaved along $Z$; it is not coherent and therefore not fine.
\item In this paper all log structures are over $k^\dagger$ even when
  this is not explicitly stated.   
  \end{enumerate}
\end{conv}

\begin{rem}
  \label{rem:joint_compatibility}
  It is a simple exercise to show that the joint compatibility
  condition is satisfied if, for all $i$, $f_i|k[u]=u^{\ell_i}$, and
  we will always assume that this is so. 
\end{rem}

\begin{dfn}
  Let $L$ be a rank-$2$ lattice and
  $S=\{(e_i, \nu_i)\mid i\in I \}$ a log datum for $L$.  A log
  structure given by wall functions $f_i$ is \emph{subordinate to $S$} if for
  all $i$
  \[
f_i=\prod_{k=1}^{k(i)} f_{i,k}
 \]
 where for all $i,k$:
 \begin{enumerate}[(i)]
 \item $Z_{i,k} =(f_{i,k}=0)\subset D_i$ is a smooth curve;
 \item $f_{i,k}|k[u]=u^{\ell_{i,k}}$.
 \end{enumerate}
\end{dfn}

For every log $S$ datum there exists a log structure that is \emph{subordinate to $S$}.

\begin{rem}
  In this paper we will always assume that the $Z_{i,k}$ are generic
  stc~(i) and~(ii). The generic attribute just means that the curves $Z_{i,k}$ are pairwise distinct and meet one another only in the origin.
\end{rem}

\subsection{Mutations of log data}
\label{sec:mutat-well-form}

In this section, we define mutations of log data. 

\begin{nota}
  For all $a\in \RR$ we write
  \[
    a_+=
    \begin{cases}
      a & \quad \text{if $a\geq 0$}\\
      0 & \quad \text{if $a<0$}
    \end{cases}
    \]
\end{nota}

\begin{dfn}
  \label{dfn:uheight}
    Let $L$ be a rank-$2$ lattice, $u\in L$ and $S=\{(e_i,\nu_i)\mid i\in
    I\}$ a log datum for $L$. The \emph{$u$-height} of $S$ is the
    natural number
    \[
h_u(S)=\sum_{i\in I}\{u,e_i\}_+.
     \]
\end{dfn}

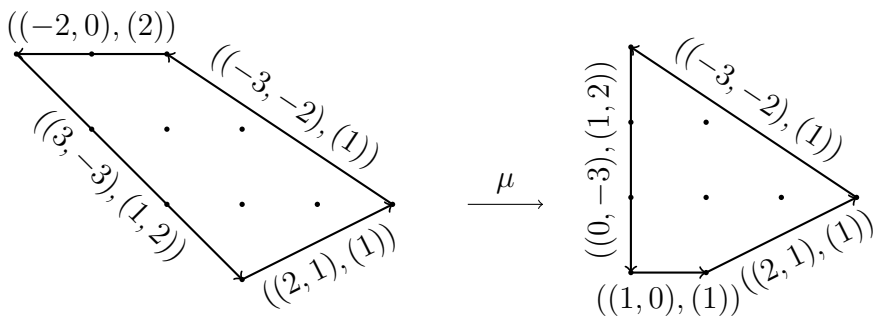
\begin{figure}[h!]
\centering
\begin{tikzpicture}
\draw[thick,->] (0,0) -- (2,1) node[midway,below,rotate=26.57]{$((2,1),(1))$};
\draw[thick,->] (2,1) -- (-1,3) node[midway,above,rotate=-33.69]{$((-3,-2),(1))$};
\draw[thick,->] (-1,3) -- (-3,3) node[midway,above,rotate=0]{$((-2,0),(2))$};
\draw[thick,->] (-3,3) -- (0,0) node[midway,below,rotate=-45]{$((3,-3),(1,2))$};
\fill (0,0) circle (1pt);
\fill (-1,1) circle (1pt);
\fill (0,1) circle (1pt);
\fill (1,1) circle (1pt);
\fill (2,1) circle (1pt);
\fill (-2,2) circle (1pt);
\fill (-1,2) circle (1pt);
\fill (0,2) circle (1pt);
\fill (-3,3) circle (1pt);
\fill (-2,3) circle (1pt);
\fill (-1,3) circle (1pt);
\draw[->] (3,1) -- (4,1) node[midway,above]{$\mut$};
\end{tikzpicture}\quad
\begin{tikzpicture}
\draw[thick,->] (0,0) -- (1,0) node[midway,below]{$((1,0),(1))$};
\draw[thick,->] (1,0) -- (3,1) node[midway,below,rotate=26.57]{$\text{ } \ \ ((2,1),(1))$};
\draw[thick,->] (3,1) -- (0,3) node[midway,above,rotate=-33.69]{$((-3,-2),(1))$};
\draw[thick,->] (0,3) -- (0,0);
\draw (-.4,1.5) node[rotate=90]{$((0,-3),(1,2))$};
\fill (0,0) circle (1pt);
\fill (1,0) circle (1pt);
\fill (0,1) circle (1pt);
\fill (1,1) circle (1pt);
\fill (2,1) circle (1pt);
\fill (3,1) circle (1pt);
\fill (0,2) circle (1pt);
\fill (1,2) circle (1pt);
\fill (0,3) circle (1pt);
\end{tikzpicture}
\caption{We represent a log datum by its polygon with each edge labelled by the corresponding $(e_i,\nu_i)$. The figure shows  a mutation in $((-2,0),(2))$.}
\label{fig:othermut}
\end{figure}

\begin{dfn}
  \label{dfn:mutation}
Let $L$ be a rank-$2$ lattice and $S=\{(e_i,\nu_i)\mid i\in I\}$ a rank-two log datum for $L$. 
We fix $j\in I$ and $1\leq k \leq k(j)$, denote $h:=h_{u_j}(S)$ and assume that $h\geq
    \ell_{j,k}$. 
    
The \emph{mutation} $\mu_{j, k}(S)$
    is the log datum for $L$ defined as follows.
To start with, let $I'$ be the same index set as $I$. We will modify $I'$ algorithmically in the following.
 \begin{enumerate}[(1)]
  \item If $i\in
    I$ and $e_i \not \in \ZZ u_j$, set $e_i^\prime = e_i + \{u_j,
    e_i\}_+u_j$ and $\nu_i^\prime = \nu_i$.
  \item
     \begin{enumerate}[(a)]
      \item If $\nu_j\neq (\ell_{j,k})$, set $e_j^\prime = (\ell_j-\ell_{j,k})u_j$ and $\nu_j^\prime = \nu_j\setminus \{\ell_{j,k}\}$. 
      \item If $\nu_j=(\ell_{j,k})$, remove $j$ from the index set $I'$.
     \end{enumerate}
  \item 
    \begin{enumerate}[(a)]
     \item
      If there exists an index $j^\star \in I$ such that
      $u_{j^\star}=-u_j$, set $e_{j^\star}^\prime = e_{j^\star}
      +(h-\ell_{j,k})u_{j^\star}$ and $\nu_{j^\star}^\prime = \nu_j\cup \{h-\ell_{j,k}\}$.
     \item If there does not exists an index $j^\star \in I$ such that
      $u_{j^\star}=-u_j$ but $d:=h-\ell_{j,k}>0$, then we add a new index $j^\star$ to $I'$ and set 
      $e'_{j^\star}=-du_j$ and $\nu'_{j^\star}=(d)$.
    \end{enumerate}
 \end{enumerate}
   We define the mutation of $S$ in $(j,k)$ to be
    \[
\mu_{j, k}(S) = \{(e_i^\prime, \nu_i^\prime) \mid i \in I'\}.
\]
\end{dfn}

Figure~\ref{fig:othermut} shows an example of a mutation that features steps (1), (2b), (3b).


\subsection{Zero-mutable log structures and statement of the 
  Conjecture}
\label{sec:0-mutable-log}

\begin{dfn}
  \label{dfn:0-mutable_logstr}
  Let $L$ be a rank-$2$ lattice.
  \begin{enumerate}[(1)]
  \item A log datum $S$ on $L$
  is \emph{zero-mutable} if there exists a sequence of
  mutations that converts $S$ into a zero-mutable rank one datum.
\item A \emph{zero-mutable log structure} is a generic log structure
  subordinate to a zero-mutable log datum. 
  \end{enumerate}
\end{dfn}

\begin{con}
  \label{con:main_conjecture}
  Let $L$ be a rank-$2$ lattice, $S$ a zero-mutable log datum,
  and $X=X_S$.
  Denote by $X^\dagger = (X, \mathfrak{M}_S)$ the log
  scheme corresponding to a generic log structure subordinate to $S$.

  Then, $X^\dagger$ admits a projective log crepant \emph{log resolution}.
  In other words there exists a log scheme $Y^\dagger$, and a projective
  morphism $\varepsilon^\dagger \colon Y^\dagger \to X^\dagger$ with the following properties:
  \begin{enumerate}[(i)]
  \item $Y=\cup_{i\in I}Y_i$ where $Y_i$ is irreducible and
    $\varepsilon$ maps $Y_i$ birationally to $X_i$.
  \item $Y^\dagger$ is log smooth over $\Spec k^\dagger$ outside finitely many points where
    $Y$ is locally \'etale isomorphic to $(xy=0) \subset
    \frac1{r}(1,-1,a,-a)$ where the log structure
    is the restriction of the divisorial log structure of this divisor embedding. Note that $Y^\dagger$ is a log smooth log
    stack at these points.
  \item The restriction of $\varepsilon$ induces the identity $Y\setminus \varepsilon^{-1}(Z) =
    X\setminus Z$ and the log structures are the same.
  \item The log canonical line bundle $K_{Y^\dagger/k^\dagger}$ is $\varepsilon$-trivial. 
  \end{enumerate}
\end{con}

\section{Examples}
\label{sec:examples}

We explain some constructions of log crepant log resolutions that give
evidence for the conjecture.

\begin{conv}
  In this section we write, for instance, $\AA^3_{x,y,z}$ to denote
  affine $3$-space $\Spec k[x,y,z]$.

  We denote by $\frac1{r}(a_1,\dots,a_n)_{x_1,\dots,x_n} $
 the quotient of $\AA^n_{x_1,\dots,x_n}$ by the $r$th roots of unity $G_r$ acting with
 weights $a_1,\dots, a_n$, that is, a generator $\zeta\in G_r$ acts by
 $\zeta \colon x_i \mapsto \zeta^{a_i} x_i$.
 
 We write the partition $1+1+1+1$ as $1^4$, etc.
\end{conv}

\subsection{What happens in codimension $1$}
\label{sec:surfaces}


\begin{lem}
  \label{lem:kink}
  Let $L$ be a rank-$2$ lattice and $S=\{(e_i, \nu_i)\mid i \in I\}$ a
log datum on $L$. Recall that, for all $i\in I$, $e_i=\ell_i u_i$ with $u_i\in L$
primitive, and $\nu_i$ is a partition of $\ell_i$. For all log
structures on $X=X_S$ compatible with the gtc structure, for all
$i\in I$, $\ell_i$ is the \emph{kink} of the log structure along $D_i\subset X$. \qed
\end{lem}

Instead of proving the lemma, which is a tautology after unpacking the
terminology, we explain what the statement means.
We consider a log structure $\fM_X$ on $X$ compatible with the gtc structure
given by a set of wall functions $\{f_i\in k[D_i] \mid i\in I\}$.

Let $\Sigma_i=\rho_i^{-1}\Sigma$ be the localisation of the fan $\Sigma$ at
$\rho_i$ and let
\[
  U_i=\Spec k[\Sigma_i] = X\setminus \Bigl(\cup_{j\neq i-1, i}
  X_j\Bigr).
\]
By choosing oriented bases $x,z=x^{e_i}$ and $z, y\in L$, we can identify
\[
U_i=(xy=0)\subset \AA^4_{x,y,z,u} \setminus(z=0) 
\]
Write $X_j^\star = X_j\cap U_i$ (for $j=i-1,i$),
$D_i^\star=D_i\cap U_i$, $f_i^\star=f_i|D_i^\star$, consider the
deformation:
\[
\pi_i \colon \fU_i=(xy-t^{\ell_i}f_i^\star=0)\subset \AA^t_{x,y,z,u,t}\setminus(z=0) \to \AA^1_t
\]
and denote by $\kappa_i\colon U_i \hookrightarrow \fU_i$ the
inclusion of the fibre $\pi_i^\star (0)$ and by $\fM_{\fU_i, U_i}$ the divisorial log structure of the log pair. Lemma~\ref{lem:kink} states that, for all $i\in I$,
\[
\fM_X|U_i = \kappa_i^\star \fM_{\fU_i, U_i}
\]

\medskip

In the next few sections we will construct several resolutions
$\varepsilon^\dagger \colon Y^\dagger \to X^\dagger$. For all these,
it will be true that, above $U_i$, the resolution is the ordinary blow
up of the curve $Z_i^\star=(f_i^\star =0)\subset D_i^\star$ on one
side only, that is, either the blow up of
$Z_i^\star \subset X_{i-1}^\star$, or the blow up of
$Z_i^\star \subset X_{i}^\star$. More precisely it will be true that
\[
  V_i=\varepsilon^{-1}(U_i) = Y_{i-1}^\star\cup Y_i^\star
\]
where
\begin{description}
\item[either] $Y_{i-1}^\star$ is the blow up of $Z_i^\star \subset
  X_{i-1}^\star$ and $Y_i^\star = X_i^\star$, 
\item[or] $Y_{i-1}^\star = X_{i-1}^\star$ and  $Y_i^\star$ is the blow up of $Z_i^\star \subset
  X_i^\star$.
\end{description}
In our context $Z_i^\star\subset D_i^\star$ is a smooth curve: it is
clear from the construction that in either case $V_i^\dagger$ is a
smooth log scheme with kink $\ell_i$ along the intersection
$Y_{i-1}^\star \cap Y_i^\star$.

\subsection{The $A_n$ singularity}
\label{sec:a_n-singularity}

\begin{dfn}
  \begin{enumerate}[(1)]
  \item Consider the vectors in $L=\ZZ^2$:
\[
  e_1=(1,0), \quad e_2= (0,n+1), \quad e_3=(-1,-n-1) 
\]
The \emph{$A_n$ log datum} is the zero-mutable log datum
\[
  S=\Bigl\{\bigl(e_1,(1)\bigr),\bigl(e_2,(1^{n+1}\bigr),\bigl(e_3,(1)\bigr)\Bigr\}.
\]
This log datum is zero-mutable by the sequence of mutations $S=\mut_{2,1}\circ\ldots\circ\mut_{2,1}\circ \mut_{1,1}(\emptyset)$ of $n+1$ many mutations in the second vector that bring the log datum down to the minimal one
$\{(e_1,(1)),(-e_1,(1))$, composed with a final mutation in the first vector.
\item The \emph{$A_n$ log structure} is the generic log structure
  subordinate to the $A_n$ log datum.
  \end{enumerate}
\end{dfn}

Let $L=\ZZ^2$ and $S$ be the $A_n$ log datum. Writing
\[
  x=x^{(1,0,0)}, \quad y=x^{(-1,-n-1,0)}, \quad z=x^{(0,1,0)}, \quad w=x^{(0,-1,0)}, \quad u=x^{(0,0,1)}
\]
in $M=L\oplus \ZZ$, we have
\[
X=
\begin{cases}
  xy-w^{n+1} &=0 \\
  zw & =0
\end{cases}
\quad \subset \AA^5_{x,y,z,w,u}
 \]
 In the notation of \S~\ref{sec:log-structures},
 \[
D_1= \AA^2_{x,u}, \quad D_2=\AA^2_{z,u}, \quad D_3=\AA^2_{y,u}
 \]
 and $X=X_1\cup X_2 \cup X_3$ where:
 \[
X_1=\AA^3_{x,z,u},\quad X_2=\AA^3_{y,z,u}, \quad X_3 = (xy-w^{n+1}=0)
\subset \AA^4_{x,y,w,u}
\]
Thus the $A_n$ log structure is given by the wall functions:
\[
f_1(x,u) = u, \quad f_2(z,u) = u^{n+1}+a_1u^nz+\cdots +a_{n+1} z^{n+1},\quad f_3(y,u)=u
\]
where $a_1,\dots, a_{n+1}$ are general constants; in particular, most
importantly, $f_2$ has $n+1$ distinct roots on $\PP^1$.

We denote by $X^\dagger$ the corresponding log scheme.

\begin{figure}[h!]
\centering
\begin{tikzpicture}[scale=2,tdplot_main_coords]
\fill (0,0,0) circle (1pt);
\fill[black!40!green,opacity=0.2] (0,0,-1.3) -- (-1.3,0,-1.3) -- (-1.3,0,1.3) -- (0,0,1.3);
\fill[black!40!red,opacity=0.2] (0,0,-1.3) -- (2.2,-1.1,-1.3) -- (2.2,-1.1,1.3) -- (0,0,1.3);
\draw[thick,black!40!green,line width=1.4pt] (0,0,0) -- (-1.3,0,1);
\draw[thick,black!40!green,line width=1.4pt] (0,0,0) -- (-1.3,0,1/3);
\draw[thick,black!40!green,line width=1.4pt] (0,0,0) -- (-1.3,0,-2/3);
\draw[thick,black!40!green,line width=1.4pt] (0,0,0) -- (-1.3,0,-1);
\draw[thick,black!40!red,line width=1.4pt] (0,0,0) -- (2.2,-1.1,0);
\draw[thick,black!40!blue,line width=1.4pt] (0,0,0) -- (2.2,1.1,0);
\fill[black!40!blue,opacity=0.2] (0,0,-1.3) -- (2.2,1.1,-1.3) -- (2.2,1.1,1.3) -- (0,0,1.3);
\draw[dotted,->] (0,0,0) -- (0,0,1.3);
\draw[dotted,->] (0,0,0) -- (-1.3,0,0);
\draw[dotted,->] (0,0,0) -- (2.2,1.1,0);
\draw[dotted,->] (0,0,0) -- (2.2,-1.1,0);
\fill (0,0,1) circle (1pt) node[left]{$u$};
\fill (2,1,0) circle (1pt) node[left]{$x$};
\fill (2,-1,0) circle (1pt) node[below]{$y$};
\fill (-1,0,0) circle (1pt) node[below]{$z$};
\fill (1,0,0) circle (1pt) node[below]{$w$};
\draw (3,2,0) node{$X_1$};
\draw (0,-1,1.3) node{$X_2$};
\draw (0,-.5,-2) node{$X_3$};
\end{tikzpicture}
\caption{Sketch of the log singular locus for the $A_3$ log structure represented inside the fan $\Sigma$, labelling the lattice points that give Stanley Reisner ring generators.}
\label{fig:An}
\end{figure}
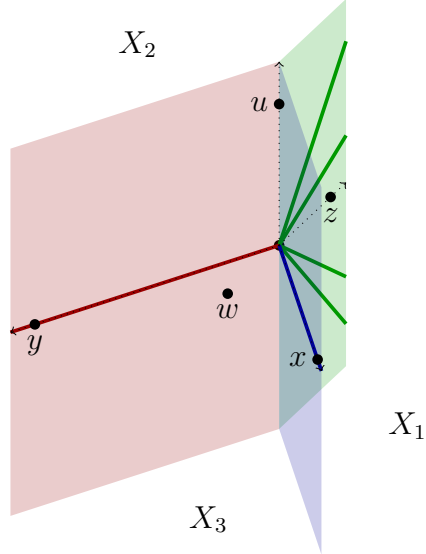

Note that, writing
\[
Z_i=(f_i=0)\subset D_i, \quad \text{and} \quad Z=Z_1\cup Z_2 \cup Z_3,
\]
we have that, strictly speaking, the wall functions $f_i$ define a log structure
$\fM_{X\setminus Z}$ on $X\setminus Z$: denoting by $i\colon
X\setminus Z \hookrightarrow X$, the log structure on $X$ is
\[
\fM_X = i_\star \fM_{X\setminus Z}
  \]

  \smallskip

Our goal in this section is to construct four log crepant log
resolutions of the $A_n$ log structure.

\begin{rem}
  The log structure on $X$ given by the wall
  functions:
  \[
    f_1(x,u) = u, \quad f_2(z,u) = u^{n+1}+az,\quad f_3(y,u)=u
  \]
   (where $a\neq 0$ is a constant) is not subordinate to the
   $A_n$ log datum, it is not zero-mutable, and it does not have a log
   crepant log resolution. This fact is closely related to the fact that
   this log structure, unlike the $A_n$ log structure, is not
   smoothable. 
 \end{rem}

The most direct way to construct log resolutions of the $A_n$ log
structure is to modify the total space of an explicit smoothing, which
is easy to write down in this case. Consider the
$4$-dimensional variety $\fX$ defined by the
equations:
\[
\fX =
\begin{cases}
  xy & = \;f_2(w,t)=w^{n+1}+a_1w^nt+\cdots +a_{n+1} t^{n+1} \\
  zw & =\;tu
\end{cases}
\quad \subset \AA^6_{x,y,z,w,u,t}
 \]
 together with the obvious morphism $\pi \colon \fX \to \AA^1_t$. It
 is clear that $X=\pi^\star (0)$ is the fibre over the origin, and it
 can be checked \cite[\S~5.4.2]{2023arXiv231213867C} that the $A_n$ log structure on $X$ is the restriction
 of the divisorial log structure $\fM_{\fX, X}$, that is, denoting by $i\colon X \hookrightarrow \fX$ the
 inclusion:
 \[
\fM_X = i^\star \fM_{\fX, X}
\]

We construct log resolutions $f^\dagger \colon Y^\dagger \to
X^\dagger$ by restricting small resolutions $\ff^\dagger
\colon \fY^\dagger \to \fX^\dagger$.

\subsubsection*{Resolutions IIIa and IIIb}
\label{sec:resol-iiia-iiib}

We construct two
log crepant log resolutions $Y^\dagger$ of $X^\dagger$; both have the property that
$Y_3$ maps isomorphically to $X_3$. These are the easiest log resolutions to construct. 

\begin{figure}[h!]
\centering
\begin{tikzpicture}[scale=2,tdplot_main_coords]
\fill (0,0,0) circle (1pt);
\fill[black!40!green,opacity=0.2] (0,0,-1.3) -- (-1.3,0,-1.3) -- (-1.3,0,1.3) -- (0,0,1.3);
\fill[black!40!red,opacity=0.2] (0,0,-1.3) -- (2.2,-1.1,-1.3) -- (2.2,-1.1,1.3) -- (0,0,1.3);
\draw[thick,black!40!green,line width=1.4pt] (-0.3,0,0) -- (-1.3,0,-.3);
\draw[thick,black!40!green,line width=1.4pt] (-0.15,0,0) -- (-1.3,0,-.7);
\draw[thick,black!40!red,line width=1.4pt] (0,0,0) -- (2.2,-1.1,0);
\draw[thick,black!40!red,line width=1.4pt,opacity=0.5] (2.2,-1.1,0) -- (1.6,-1.1,0) -- (-0.6,0,0);
\fill[black!40!red,opacity=0.4] (2.2,-1.1,0) -- (1.6,-1.1,0) -- (-0.6,0,0) -- (0,0,0);
\draw[thick,black!40!blue,line width=1.4pt] (0,0,0) -- (2.2,1.1,0);
\draw[thick,black!40!blue,line width=1.4pt,opacity=0.5] (2.2,1.1,0) -- (1.6,1.1,0) -- (-0.6,0,0) -- (0,0,0);
\fill[thick,black!40!blue,opacity=0.3] (2.2,1.1,0) -- (1.6,1.1,0) -- (-0.6,0,0) -- (0,0,0);
\draw[thick,black!40!green,line width=1.4pt] (-0.45,0,0) -- (-1.3,0,.7);
\draw[thick,black!40!green,line width=1.4pt] (-0.6,0,0) -- (-1.3,0,.3);
\fill[black!40!blue,opacity=0.2] (0,0,-1.3) -- (2.2,1.1,-1.3) -- (2.2,1.1,1.3) -- (0,0,1.3);
\draw[dotted,->] (0,0,0) -- (0,0,1.3);
\draw[dotted,->] (0,0,0) -- (-1.3,0,0);
\draw[dotted,->] (0,0,0) -- (2.2,1.1,0);
\draw[dotted,->] (0,0,0) -- (2.2,-1.1,0);
\fill (0,0,1) circle (1pt) node[left]{$u$};
\fill (2,1,0) circle (1pt) node[left]{$x$};
\fill (2,-1,0) circle (1pt) node[below]{$y$};
\fill (-1,0,0) circle (1pt) node[below]{$z$};
\fill (1,0,0) circle (1pt) node[below]{$w$};
\draw (3,2,0) node{$X_1$};
\draw (0,-1,1.3) node{$X_2$};
\draw (0,-.5,-2) node{$X_3$};
\end{tikzpicture}
\caption{Sketch of the central fiber $t=0$ inside the intermediate space $\fY'$ in the IIIa resolution of the $A_3$ log structure. 
The inverse image of $Z_1$ and $Z_3$ gives a $\PP^1$-bundle respectively while
the strict transform of $Z_2$ is a disjoint union of four lines, so a smooth curve.}
\label{fig:An-res}
\end{figure}
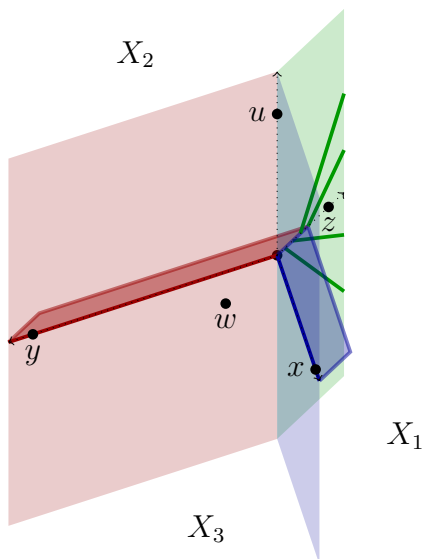

Start off with the blow up $\varepsilon \colon \fY' \to
\fX$ of the ideal $(z,u)$, so $\fY'=(u\zeta -z \nu=0)\subset \fX \times
\PP^1_{[\nu:\zeta]}$. The chart most relevant for us is the chart
where $\nu = 1$: in this chart, we have $z=u\zeta$ and
hence, solving for $z$ and for $t=\zeta w$ we get
\[
\fY' = \bigl(xy=w^{n+1}f_2(1, \zeta)\bigr) \subset \AA^5_{x,y,\zeta,w,u}.
\]
A sketch is shown in Figure~\ref{fig:An-res}.
The divisorial log structure $\fM_{\fY', Y}$ is log smooth outside the
smooth and possibly reducible curve
\[(f_2(1,\zeta)=0) \subset (Y_1\cap Y_2)\setminus Y_3
\]
and a small log resolution $\fY\to \fY'$ can then be obtained by blowing up this curve in $Y_2$
only (IIIa) or in $Y_1$ only (IIIb).

\subsubsection*{Resolution IIc}
\label{sec:resolution-iic}

We construct a log crepant log resolution $Y^\dagger$ of $X^\dagger$,
which we call IIc, with the property that $Y_2$ is mapped
isomorphically to $X_2$.

First off let
$\varepsilon \colon \FF \to \AA^6_{x,y,z,w,u,t}$ be the weighted blow up with
weights $0,n+1,0,1,0,1$. The chart of $\FF$ most relevant for us is isomorphic
to 
$
\frac1{n+1}(0,1,0,-1,0,-1)_{x, \eta, z, w^\prime,u,t^\prime}
$
and the morphism $\varepsilon$ is given by
\[
(x,y,z,w,u,t) \mapsto (x, \eta^{n+1},z, w^\prime \eta, u, t^\prime \eta),
\]
hence
\[
  \fY
  =
  \begin{cases}
    x  & = \; f_2(w^\prime, t^\prime) \\
    zw^\prime & = \;  t^\prime u 
  \end{cases}
  \qquad
  \subset\qquad
  \frac1{n+1}(0,1,0,-1,0,-1)_{x, \eta, z, w^\prime,u,t^\prime}
  \]
  in this chart. Solving for $x$, we simplify to
  \[
\fY = (zw^\prime - t^\prime u=0) \quad\subset\quad
\frac1{n+1}(1,0,-1,0,-1)_{\eta, z, w^\prime,u,t^\prime} 
\]
and the morphism to $\AA^1_t$ is given by $t=t^\prime \eta$.
The fibre $Y$ over $0\in \AA^1_t$ thus consist of three components:
$Y=Y_1\cup Y_2 \cup Y_3$
where
\begin{multline*}
  Y_1=(zw^\prime-t^\prime u=0) \subset
  \frac1{n+1}(0,-1,0,-1)_{z,w^\prime,u, t^\prime} \;,\\
  Y_2=\frac1{n+1}(1,0,0)_{\eta, z, u}\;,\quad Y_3=\frac1{n+1}(1,-1,0)_{\eta, w^\prime,u}
\end{multline*}
and $\varepsilon$ maps $Y_2$ isomorphically to $X_2$, it maps $Y_3$
isomorphically to $X_3$, and it maps $Y_1$ birationally to $X_1$. The
morphism $\varepsilon_{|Y_1} \colon Y_1 \to X_1=\AA^3_{x,z,u}$ is
given by $x\mapsto f_2(w^\prime, t^\prime)$. One can check that, over
the generic point of the curve $Z_1\cup Z_2$, the
morphism $\varepsilon_{|Y_1}$ is isomorphic to the (ordinary) blow-up
of $Z_1\cup Z_2\subset X_1$; however, something interesting is going on over
$0\in X_1$.

The log structure $\fM_{\fY, Y}$ is log smooth away from the
$\eta$-axis; to resolve it, we need a small modification of $\fY$ that
blows up the $\eta$-axis. There are two ways of doing this: blowing up
the locus $z=u=0$, or blowing up the locus $w^\prime = u = 0$. The
first of these, it can be seen, leads back to IIIa.

Hence consider the (toric) blow up
\[
  \widetilde{\varepsilon} \colon \widetilde{\FF} \to
  \frac1{n+1}(1,0,-1,0,-1)_{\eta, z, w^\prime, u, t^\prime} \subset \FF
\]
of the ambient space along the locus $w^\prime = u = 0$.

The first chart for $\widetilde{\FF}$ is isomorphic to
$\frac1{n+1}(1,0,-1,1,-1)_{\eta,z,\omega,u^\prime,t^\prime}$
and the morphism $\widetilde{\varepsilon}$ is given by
\[
(\eta, z, w^\prime,u,t^\prime) \mapsto (\eta, z, \omega,u^\prime \omega, t^\prime)
\]
In this chart, the proper transform $\widetilde{\fY}$ is given by the
equation $z=t^\prime u^\prime$; solving for $z$ we see that
\[
  \widetilde{\fY}=\frac1{n+1}(1,-1,1,-1)_{\eta,\omega,u^\prime,t^\prime}
  \quad \text{with} \quad t=t^\prime \eta \colon \widetilde{\fY} \to \AA^1_t.
\]
We see from this description that the log structure
$\fM_{\widetilde{\fY},\widetilde{Y}}$ is not log smooth at the origin: it has a
singularity as in Conjecture~\ref{con:main_conjecture}(ii).

The second chart for $\widetilde{\FF}$ is isomorphic to
$\frac1{n+1}(1,0,-1,0,-1)_{\eta,z,\widetilde{w},\nu,t^\prime}$
and the morphism $\widetilde{\varepsilon}$ is given by
\[
(\eta, z, w^\prime,u,t^\prime) \mapsto (\eta, z, \widetilde{w}\nu,\nu, t^\prime)
\]
In this chart, the proper transform $\widetilde{\fY}$ is given by the
equation $t^\prime =z\widetilde{w}$; solving for $t^\prime$ we see that
\[
  \widetilde{\fY}=\frac1{n+1}(1,0,-1,0)_{\eta,z,\widetilde{w},\nu}
  \quad \text{with} \quad t=\eta z\widetilde{w} \colon \widetilde{\fY} \to \AA^1_t.
\]
We see from this description that the log structure
$\fM_{\widetilde{\fY},\widetilde{Y}}$ is log smooth in this chart.

\subsubsection*{Resolution Ic}
\label{sec:resolution-iic}

There is a similar construction of a log crepant log resolution $Y^\dagger$ of $X^\dagger$,
which we call Ic, with the property that $Y_1$ is mapped
isomorphically to $X_1$. It can be obtained by post-composing the resolution IIc with automorphism of $X^\dagger$ that interchanges $x$ and $y$ while leaving the other variables fixed.

\subsection{Tom and Jerry}
\label{sec:tom-jerry}

A.C.\ heard Miles Reid talk about Tom and Jerry in the late 90s as two
constructions of codimension-$4$ Gorenstein rings that are frequently met
in practice; Miles suggests that they may have appeared in print for
the first time in~\cite{MR2956040}.

\subsubsection*{Tom and Jerry log structures}
\label{sec:tom-jerry-log}

\begin{dfn}
  Let $L=\ZZ^2$ and consider the vectors
\[
e_1=(3,0),\quad e_2=(0,2),\quad e_3=(-3,-2).
\]
  \begin{enumerate}[(1)]
  \item The \emph{Tom log datum} is the log datum
    \[ S=\Bigl\{ \bigl(e_1, (2,1)\bigr), \bigl(e_2, (1,1)\bigr), \bigl(e_3,(1)\bigr)\Bigr\}.
    \]
    A \emph{Tom log structure} is a generic log structure
    subordinate to the Tom log datum.
  \item The \emph{Jerry log datum} is the log datum
    \[ S=\Bigl\{\bigl(e_1, (1,1,1)\bigr), \bigl(e_2, (2)\bigr), \bigl(e_3,(1)\bigr)\Bigr\}.
    \]
    A \emph{Jerry log structure} is a generic log structure
    subordinate to the Jerry log datum.
  \end{enumerate}
\end{dfn}

\begin{figure}[h!]
\centering
\begin{tikzpicture}
\draw[thick,->] (0,0) -- (3,0) node[midway,below]{$((3,0),(\textbf{2},1))$};
\draw[thick,->] (3,0) -- (3,2);
\draw (3.4,1) node[rotate=-90]{$((0,2),(1,1))$};
\draw[thick,->] (3,2) -- (0,0) node[midway,above,rotate=33.69]{$((-3,-2),(1))$};
\draw[->] (4,1) -- (5,1);
\end{tikzpicture}
\begin{tikzpicture}
\draw[thick,->] (0,0) -- (1,0) node[midway,below]{$((1,0),(1))$};
\draw[thick,->] (1,0) -- (3,2) node[midway,below,rotate=45]{$((2,2),(1,1))$};
\draw[thick,->] (3,2) -- (0,0) node[midway,above,rotate=33.69]{$((-3,-2),(1))$};
\draw (4,1) node{\large$\simeq$};
\end{tikzpicture}
\begin{tikzpicture}
\draw[thick,->] (0,0) -- (1,0) node[midway,below]{$((1,0),(1))$};
\draw[thick,->] (1,0) -- (1,2);
\draw (1.4,1) node[rotate=-90]{$((0,2),(1,1))$};
\draw[thick,->] (1,2) -- (0,0) node[midway,above,rotate=63.43]{$((-1,-2),(1))$};
\end{tikzpicture}
\caption{Tom log datum and mutation to the $A_1$ log datum with $\ell_{j,k}$ in bold.}
\label{fig:tommut}
\end{figure}
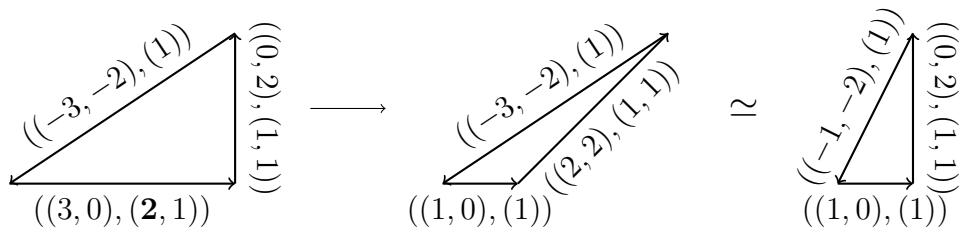

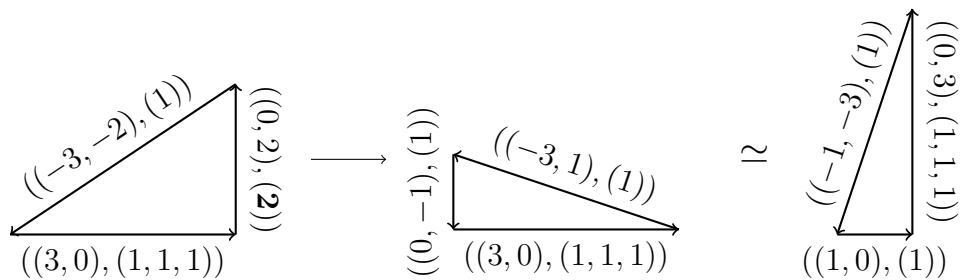
\begin{figure}[h!]
\centering
\begin{tikzpicture}
\draw[thick,->] (0,0) -- (3,0) node[midway,below]{$((3,0),(1,1,1))$};
\draw[thick,->] (3,0) -- (3,2);
\draw (3.4,1) node[rotate=-90]{$((0,2),(\textbf{2}))$};
\draw[thick,->] (3,2) -- (0,0) node[midway,above,rotate=33.69]{$((-3,-2),(1))$};
\draw[->] (4,1) -- (5,1);
\end{tikzpicture}
\begin{tikzpicture}
\draw[thick,->] (0,0) -- (3,0) node[midway,below]{$((3,0),(1,1,1))$};
\draw[thick,->] (3,0) -- (0,1) node[midway,above,rotate=-18.43]{$((-3,1),(1))$};
\draw[thick,->] (0,1) -- (0,0);
\draw (-.4,.5) node[rotate=90]{$((0,-1),(1))$};
\draw (4,1) node{\large$\simeq$};
\end{tikzpicture}
\begin{tikzpicture}
\draw[thick,->] (0,0) -- (1,0) node[midway,below]{$((1,0),(1))$};
\draw[thick,->] (1,0) -- (1,3);
\draw (1.4,1.5) node[rotate=-90]{$((0,3),(1,1,1))$};
\draw[thick,->] (1,3) -- (0,0) node[midway,above,rotate=71.57]{$((-1,-3),(1))$};
\end{tikzpicture}
\caption{Jerry log datum and mutation to the $A_2$ log datum with $\ell_{j,k}$ in bold.}
\label{fig:jerrymut}
\end{figure}

The reader can check that the Tom and the Jerry log data
are zero-mutable; Figures \ref{fig:tommut},\,\ref{fig:jerrymut} illustrate a possible first step in a sequence of mutations to the empty log datum. 

The underlying scheme $X=X_S$ is the same for the Tom and Jerry log
data: $X=X_1\cup X_2 \cup X_3$, and
\[
X_1= \AA^3, \quad X_2=\frac1{3}(1,1,0),\quad X_3=\frac1{2}(1,1,0)
\]
and hence the embedding codimension is $4$. 

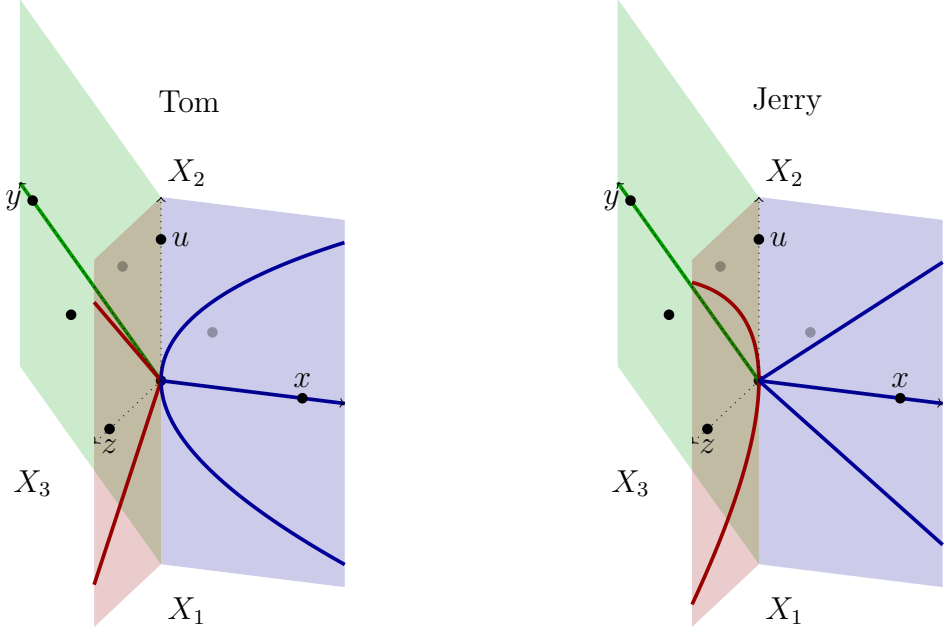
\begin{figure}[h!]
\centering
\begin{tikzpicture}[scale=2,tdplot_main_coords]
\draw (0,.2,2) node{Tom};
\fill (0,0,0) circle (1pt);
\fill[black!40!green,opacity=0.2] (0,0,-1.3) -- (-3.3,-2.2,-1.3) -- (-3.3,-2.2,1.3) -- (0,0,1.3);
\fill[black!40!blue,opacity=0.2] (0,0,-1.3) -- (0,1.3,-1.3) -- (0,1.3,1.3) -- (0,0,1.3);
\draw[thick,black!40!green,line width=1.4pt] (0,0,0) -- (-3.3,-2.2,0);
\draw[thick,black!40!red,line width=1.4pt] (0,0,0) -- (1.3,0,1);
\draw[thick,black!40!red,line width=1.4pt] (0,0,0) -- (1.3,0,-1);
\draw[thick,black!40!blue,line width=1.4pt,smooth,samples=100,domain=0:1.14] plot (0,\x^2,\x);
\draw[thick,black!40!blue,line width=1.4pt,smooth,samples=100,domain=-1.14:0] plot (0,-\x^2,\x);
\draw[thick,black!40!blue,line width=1.4pt] (0,0,0) -- (0,1.3,0);
\fill[black!40!red,opacity=0.2] (0,0,-1.3) -- (1.3,0,-1.3) -- (1.3,0,1.3) -- (0,0,1.3);
\draw[dotted,->] (0,0,0) -- (0,0,1.3);
\draw[dotted,->] (0,0,0) -- (-3.3,-2.2,0);
\draw[dotted,->] (0,0,0) -- (1.3,0,0);
\draw[dotted,->] (0,0,0) -- (0,1.3,0);
\fill (-1,-1,0) circle (1pt);
\fill[opacity=.3] (-1,0,0) circle (1pt);
\fill[opacity=.3] (-2,-1,0) circle (1pt);
\fill (0,0,1) circle (1pt) node[right]{$u$};
\fill (-3,-2,0) circle (1pt) node[left]{$y$};
\fill (1,0,0) circle (1pt) node[below]{$z$};
\fill (0,1,0) circle (1pt) node[above]{$x$};
\draw (-.5,0,-1.8) node{$X_1$};
\draw (-.5,0,1.3) node{$X_2$};
\draw (-3,-2,-2) node{$X_3$};
\end{tikzpicture}
\hspace{3cm}
\begin{tikzpicture}[scale=2,tdplot_main_coords]
\draw (0,.2,2) node{Jerry};
\fill (0,0,0) circle (1pt);
\fill[black!40!green,opacity=0.2] (0,0,-1.3) -- (-3.3,-2.2,-1.3) -- (-3.3,-2.2,1.3) -- (0,0,1.3);
\fill[black!40!blue,opacity=0.2] (0,0,-1.3) -- (0,1.3,-1.3) -- (0,1.3,1.3) -- (0,0,1.3);
\draw[thick,black!40!green,line width=1.4pt] (0,0,0) -- (-3.3,-2.2,0);
\draw[thick,black!40!red,line width=1.4pt,smooth,samples=100,domain=0:1.14] plot (\x^2,0,\x);
\draw[thick,black!40!red,line width=1.4pt,smooth,samples=100,domain=-1.14:0] plot (-\x^2,0,\x);
\draw[thick,black!40!blue,line width=1.4pt] (0,0,0) -- (0,1.3,1);
\draw[thick,black!40!blue,line width=1.4pt] (0,0,0) -- (0,1.3,0);
\draw[thick,black!40!blue,line width=1.4pt] (0,0,0) -- (0,1.3,-1);
\fill[black!40!red,opacity=0.2] (0,0,-1.3) -- (1.3,0,-1.3) -- (1.3,0,1.3) -- (0,0,1.3);
\draw[dotted,->] (0,0,0) -- (0,0,1.3);
\draw[dotted,->] (0,0,0) -- (-3.3,-2.2,0);
\draw[dotted,->] (0,0,0) -- (1.3,0,0);
\draw[dotted,->] (0,0,0) -- (0,1.3,0);
\fill (-1,-1,0) circle (1pt);
\fill[opacity=.3] (-1,0,0) circle (1pt);
\fill[opacity=.3] (-2,-1,0) circle (1pt);
\fill (0,0,1) circle (1pt) node[right]{$u$};
\fill (-3,-2,0) circle (1pt) node[left]{$y$};
\fill (1,0,0) circle (1pt) node[below]{$z$};
\fill (0,1,0) circle (1pt) node[above]{$x$};
\draw (-.5,0,-1.8) node{$X_1$};
\draw (-.5,0,1.3) node{$X_2$};
\draw (-3,-2,-2) node{$X_3$};
\end{tikzpicture}
\caption{Log singular loci for Tom and Jerry.}
\label{fig:tomjerry}
\end{figure}

In this section, we sketch constructions of log crepant log
resolutions of the Tom and Jerry log structures. Denoting by
$X^\dagger$ either the Tom or the Jerry log structure, an explicit
smoothing $f \colon \fX \to \AA^1$ of $X^\dagger$ can be
constructed~\cite{MR4307202}.\footnote{Log structures do not feature
  in that paper. That paper constructs two smoothing families of $X$,
  which it calls the Tom family and the Jerry family. It can be
  checked that the restrictions to $X$ of the divisorial log
  structures of the two families are Tom and Jerry log structures,
  which also explains our terminology.} In
principle it should be possible to resolve $X^\dagger$ by a small
modification of the total space $\fX$. However, as detailed
in~\cite{MR4307202}, the affine coordinate ring of $\fX$ is a 
codimension $4$ Gorenstein ring defined by 
 $9$ equations and $16$
syzygies. This structure
 is typical for codimension-$4$ Gorenstein
rings that one frequently meets in practice. 
Constructing small modifications of $\fX$ through explicit
blow-ups of its ambient space  and analyzing the results, as we did for the $A_n$
singularity, is not the most practical way for obtaining log resolutions of  $X^\dagger$. 
Instead, it is better to abandon any reference to the total space and directly modify the log scheme $X^\dagger$.

\subsubsection*{Resolving the Tom log structures}
\label{sec:resolving-tom-log}

The key observation is that if $S$ is the Tom log datum, then
\[
\nu_{1,1} (S) = \overline{S}
\]
is the $A_1$ log datum, see Figure~\ref{fig:tommut}. It follows that $X_1=\overline{X}_1$ and $X_3\cong
\overline{X}_3$. Let $\overline{Y}^\dagger$ be the log resolution
IIc. We obtain a resolution of $X^\dagger$ by gluing
$\overline{Y}_1\cup \overline{Y}_3$ to $X_2$ in the obvious way.

\subsubsection*{Resolving the Jerry log structures}
\label{sec:resolving-jerry-log}

Let $S$ be the Jerry log datum. It is true that
\[
\nu_{2,1} (S) = S^\prime
\]  
is (isomorphic to) the $A_2$ log datum, as shown in Figure~\ref{fig:jerrymut}, but this is not super-helpful
for constructing a log resolution of the Jerry log structures. Here we
sketch the construction of a resolution, leaving the details to the
reader. 

Writing
\[
  x=x^{(1,0,0)},\quad z=x^{(0,1,0)}, \quad y=x^{(-3,-2,0)}, \quad
  u=x^{(0,0,1)}
\]
in $M=L\oplus \ZZ$, we have
\[
  D_1 = \AA^2_{x,u},\quad D_2 = \AA^2_{z,u},\quad D_3=\AA^2_{y,u}
\]
The Jerry log structure is given by the wall functions:
\[
f_1=u^3 +a_1u^2x+a_2ux^2+ a_3x^3, \quad f_2=u^2+az, \, f_3=u 
\]
where $a, a_1, a_2, a_3$ are generic constants. As customary, we
let $Z_i=(f_i=0)\subset D_i$.

As a first step, we let:
\begin{enumerate}[(1)]
\item $Y_1\to X_1$ be the blow up of the curve
$Z_2\subset X_1$;
\item $Y_2=X_2$;
\item $Y_3 \to X_3 = \frac1{2}(1,1,0)$ the (toric)
blow up of the curve $Z_3\subset X_3$. 
\end{enumerate}

Glue $Y_1$, $Y_2$ and $Y_3$ in the obvious way to get $Y=Y_1\cup Y_2
\cup Y_3$; this is a gtc space and the results
of~\cite{2023arXiv231213867C} allow to put a log structure on $Y$ and
to promote the birational mophism $\varepsilon \colon Y \to X$ to a
morphism of log schemes $\varepsilon^\dagger \colon Y^\dagger \to
X^\dagger$.

The log scheme $Y^\dagger$ is singular along the proper transform
$Z_1^\prime \subset D_1^\prime =Y_1\cap Y_3 \subset Y$. There is a
unique ordinary quadratic singular point $q\in D_1^\prime$, and
$q\in Z_1^\prime\subset Y_1$ looks like three coordinate axes meeting at $0\in \AA^3$. Let
$\widetilde{Y}_1 \to Y_1$ be the extremal blow up of $Z_1^\prime
\subset Y_1$. The resolution we want is obtained gluing
$\widetilde{Y}_1$ back to $Y_2\cup Y_3$; the resulting log scheme has a
$\frac1{2}(1,1,1,1)$ singularity as in Conjecture~\ref{con:main_conjecture}(ii).

\section{Relation to~\cite{MR4381899} and the Gross--Siebert program}
\label{sec:relation-}

We sketch 
connections with the work of Corti--Filip--Petracci~\cite{MR4381899}
and the Gross--Siebert program~\cite{MR2846484, MR3415066, MR4520304, MR4462625}.

Let $L^\star$ be a rank-$2$ lattice and $F\subset L^\star$ an integral
polygon. Consider the cone over F at height one:
\[
\sigma_F = \langle F \times \{1\}\rangle_+ \subset N=L^\star \oplus \ZZ 
\]
The central topic of~\cite{MR4381899} is the deformation theory of the
toric affine \mbox{$3$-fold} $X_{\sigma_F}$. The main point of that
paper is to conjecture that there is a bijective correspondence
between smoothing components of $X_{\sigma_F}$ and zero-mutable Laurent
polynomials with Newton polygon $F$.

For the present discussion it is not necessary to recall the
definition of zero-mutable Laurent polynomials. Let $u_i\in L$ be the
inner unit normals of $F$ and $\ell_i$ the integral length of the
corresponding edge and set $e_i=\ell_iu_i$. It can easily be shown that
there is a bijective correspondence between zero-mutable Laurent
polynomials with Newton polytope $F$ and zero-mutable log data
$S=\{(e_i, \nu_i)\}$. Note that $X_S$ is a Mumford degeneration of
$X_{\sigma_F}$. Now choose a zero-mutable log datum $S$ and zero-mutable log
structure $\fM_S$ and let $X^\dagger = (X,\fM_S)$. In this paper we conjecture the existence of log
resolutions $\varepsilon^\dagger \colon Y^\dagger \to X^\dagger$. We
expect that $Y^\dagger$ has smoothing deformations which in turn induce smoothing deformations of $X^\dagger$, and that these smoothing
deformations of $X^\dagger$ arise naturally as Mumford degenerations of the smoothing
deformations of $X_{\sigma_F}$ in the smoothing component associated,
via the conjecture of~\cite{MR4381899}, to the zero-mutable
Laurent polynomial corresponding to the zero-mutable log datum.

All this ties in well with the Gross--Siebert program. Starting from a
zero-mutable log datum $S=\{(e_i,\nu_i)\}$ and subordinate log
structure given by a set of functions $\{f_i\}$, construct an
initial scattering diagram on $M=L\oplus \ZZ$ with initial walls
$\rho_i$ by placing the function $f_i$ on the wall $\rho_i$. We expect
to be able to run the Kontsevich--Soibelman/Gross--Siebert program on
this scattering diagram (at the moment this does not follow from
any results in the literature), and that this will produce a
distinguished $1$-parameter (formal) smoothing of $X^\dagger$. This,
we expect, can be simultaneously log resolved, inducing the log
resolution $\varepsilon^\dagger \colon Y^\dagger \to X^\dagger$ whose
existence is conjectured in this paper.

\subsection{Higher dimensions}
We expect zero-mutable log structures can be defined and studied in any dimension and be subject to a conjecture similar to Conjecture~\ref{con:main_conjecture} though we do not know at this point what exceptional points or loci we should expect in the crepant resolution, e.g., whether we can generally expect a log smooth stack - most likely not. There are two steps of generalization.
\subsubsection{Step I} The lattice $L$ remains rank two but $k[u]$ gets replaced by $k[u_1,...,u_r]$. The log datum $S$ then decomposes disjointly $S=S_1\cup...\cup S_r$ and each $S_i$ is itself a log datum literally as defined before except with $u$ replaced by $u_i$ for $i=1,...,r$. Potentially, we also want to permit an abelian quotient under a finite group that acts by multiplication with roots of unity on the $u$-coordinates.
\subsubsection{Step II} If the rank of the lattice $L$ increases, we need to additionally specify the fan $\Sigma$ in $L$. It was uniquely determined by its rays when the rank of $L$ is two but there is more information in it when the rank of $L$ goes up.
 The fan $\Sigma$ shall permit a piecewise linear strictly convex function, or equivalently, the fan shall be the normal fan of a dual Newton polytope. This additional requirement corresponds to the Gorenstein condition. Each two-face of the Newton polytope is dual to a cone $\tau\in \Sigma$ so that the localization quotient $\tau^{-1}\Sigma/\tau$ lives in a rank two quotient of $L$. With this quotient, we are now back to the situation of Step I. However, when $\tau$ is larger than dimension two, its dual gives a new situation that we have not studied yet.

\bibliographystyle{plainurl}
\bibliography{LogRes}

\end{document}